\documentclass{article}

\usepackage[english]{babel}
\usepackage[letterpaper,top=0.74in,bottom=2cm,left=0.5in,right=0.5in,marginparwidth=0.5in]{geometry}

\makeatletter
\if@titlepage
  \renewenvironment{abstract}{%
      \titlepage
      \null\vfil
      \@beginparpenalty\@lowpenalty
        \paragraph{\abstractname:}
        \@endparpenalty\@M
      }%
     {\par\vfil\null\endtitlepage}
\else
  \renewenvironment{abstract}{%
      \if@twocolumn
        \section*{\abstractname}%
      \else
        \small
        \paragraph{\abstractname:}
      \fi}
      {\if@twocolumn\else\par\bigskip\fi}
\fi
\makeatother

\usepackage{amsmath}
\usepackage{graphicx}
\usepackage[colorlinks=true, allcolors=blue]{hyperref}

\title{Introducing the Role of Shaping Order K in Set Shaping Theory}
\author{Sochima Biereagu\footnote{Author: sochima.eb@gmail.com}}
\date{}

\begin{document}
\maketitle

\abstract{
Set Shaping Theory, an emerging area of study, delves into the transformation of data sets via bijection functions \cite{set_shaping_theory}. Central to this theory is the parameter $K$, which determines the extent of transformation, essentially reshaping the data. This article introduces the pivotal role of $K$ in the Set Shaping Theory, shedding light on its implications for data compression and transformation dynamics.
}

\section*{Introduction}
Within the Set Shaping Theory, the parameter $K$ is a pivotal determinant of data transformation. A profound analytical understanding of $K$ is crucial to harness the full potential of the theory. In this section, we embark on an introductory exploration of $K$ and its implications

\subsubsection*{Brief Overview of Set Shaping Theory}
Fundamentally, the Set Shaping Theory investigates bijection functions that transform one set of strings into another. The crux is the transformation of string sets, such as $X^N$, into longer string sets, like $Y^{N+K}$, with each string in the latter being $K$ symbols longer than its counterpart in the former.

\section*{The Mathematical Foundation of $K$}

\subsection*{Definitions and Basic Concepts}
\begin{itemize}
    \item \textbf{Source Representation}: A source, denoted as $X$, is characterized by the triplet $X = (x; A; P)$, where:
    \begin{itemize}
        \item $x$ is the value of a random variable.
        \item $A$ represents the possible states or values of $x$.
        \item $P$ is the probability distribution associated with the states.
        \item We call $X^N$ the set of all possible strings $x = x_1, \ldots, x_j, \ldots, x_N$ generated by $X$.
    \end{itemize}
    \item \textbf{Entropy \cite{Shannon}}: A measure of unpredictability or randomness of a set, $H(X)$, is given by: 
    $$H(X) = -\sum_i p_i \log_b p_i$$
    \item \textbf{Information Content (empirical entropy)}: Provides a measure of the "weight" or "significance" of a specific string within the set, For a specific string $x_i$ in $X^N$, its information content, $I(x_i)$, is:
    $$I(x_i) = -\sum_{j=1}^N \log_2 p(x_j)$$
\end{itemize}

The bijection function $f$ that performs the transform is defined as:
$$f: X^N \rightarrow Y^{N+K}$$
\newline
\textbf{Transformation process:}
\begin{itemize}
    \item Using the bijection function $f$, a string $x$ is mapped to a new string $y$ in $Y^{N+K}$
    \item The key insight from set shaping theory is that the strings in $Y^{N+K}$ are chosen based on their lower information content from the set $X^{N+K}$. This means that, even though the strings in $Y^{N+K}$ are longer, they might be more "compressible" due to their reduced information content.
\end{itemize}

\begin{figure}[h]
    \centering
    \includegraphics[width=0.5\linewidth]{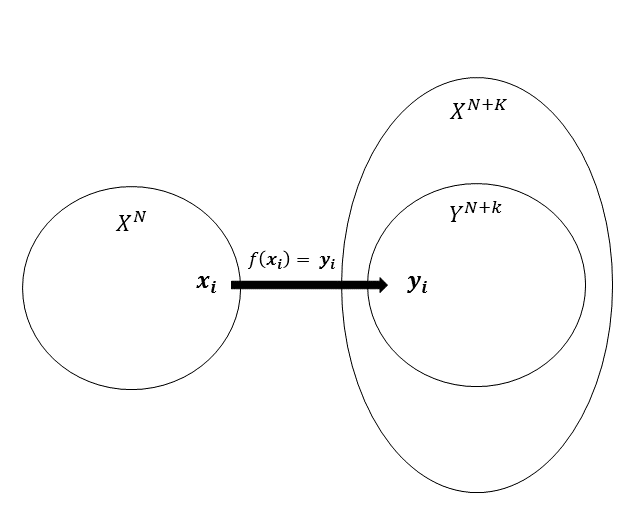}
    \caption{Set shaping theory bijection illustration $f(x_i)$$^{[3]}$}
    \label{fig:1}
\end{figure}

\subsection*{The Role of Shaping Order $K$}

The shaping order, $K$, is central to the transformation dynamics in the Set Shaping Theory.
It acts as a determinant of how much length is added to each string during the transformation. If you're transforming a string from the set $X^N$ using the function $f$, the resulting string in $Y^{N+K}$ will be $K$ characters/symbols longer.
\newline\newline
Additionally, the magnitude of $K$ can be seen as an indicator of the depth or extent of the transformation. A larger $K$ suggests a more profound reshaping of the data, potentially leading to more strings with even less information content.
\newline
This balance between the added length and the reduction in information content is crucial for the practical application of the Set Shaping Theory, especially in data compression scenarios. 

\subsection*{Case Study}

To truly grasp the significance and application of the $K$ parameter in Set Shaping Theory (SST), consider the following example
\newline
Given a source defined by an ensemble $X=(x; A; P)$ with a uniform probability distribution, we'll show the results of applying the SST with different values of $K$ from $0$ to $7$.
\newline
Let $A = \{0,1,2\}\ and\ N = 10$. This implies that $X^N$ would contain all possible strings of length $10$ produced by the source $X$. That would be $3^{10} = 59,049$ strings.
\subsubsection*{Some Definitions:}
\begin{itemize}
\item The probability $P(x_i)$ that the source $X$ generates the sequence $x_i$ is given as:
$$P(x_i) = \prod_{j=1}^N p(x_j)$$

\item We call the average information content of a sequence generated by a source $X=(x; A; P)$ the summation of the product between the information content of the sequences belonging to $X^N$ is their probability:
$$I(x) = \sum_{i=1}^{|X|^N} P(x_i)I(x_i)$$

\item Because our bijection function $f$ transforms the strings $x \in X^N$ into the strings $y \in Y^{N+K}$ consequently, the average information content changes as follows.
$$I(y) = \sum_{i=1}^{|X|^N} P(x_i)I(y_i)$$
\end{itemize}

For this test, we'll compare the average information content of the set $Y^{N+K}$ for the values of $K$ from $0$ to $7$. When $K = 0$ it's the same as not applying the SST, because the string length doesn't increase.
\newline
The first column shows the value of $K$, the second column shows the length of the strings in the set $Y^{N+K}$ and the last column shows the average information content of the strings in set $Y^{N+K}$.

\begin{table}[ht]
    \centering
    \begin{tabular}{|c|c|c|}
        \hline
        K & \text{Length of strings in $Y^{N+K}$} & \text{$I(y)$} \\
        \hline
        0 & 10 & 14.263 \\
        \hline
        1 & 11 & 14.136 \\
        \hline
        2 & 12 & 14.006 \\
        \hline
        3 & 13 & 13.694 \\
        \hline
        4 & 14 & 13.322 \\
        \hline
        5 & 15 & 13.612 \\
        \hline
        6 & 16 & 13.809 \\
        \hline
        7 & 17 & 13.969 \\
        \hline
    \end{tabular}
    \caption{\textit{The average information content $I(y)$ in bits calculated for strings in the set $Y^{N+K}$, \text{for varying values of } $K$}}
    \label{tab:my_label}
\end{table}
 
In this table, we can observe that as $K$ increases, even though the length of the strings in $Y^{N+K}$ increases, the average information content decreases. This showcases the essence of the Set Shaping Theory and the role of $K$ in modulating the information content of transformed sets.

\section*{Implications and Applications}

\subsection*{Impact on Transformation Dynamics}

The parameter $K$ plays a pivotal role in determining the nature and outcome of the transformation dynamics in SST:

\begin{itemize}
    \item \textbf{Influence on String Length:} As highlighted in our example, $K$ directly determines the added length to each string in the source $X$. A larger $K$ results in longer strings in the transformed set $Y$.
    \item \textbf{Alteration of Information Content:} While the addition of length might seem counter-intuitive, the crux lies in the reshaped information content. SST ensures that the strings in $Y$ have reduced information content, making them potentially more compressible.
    \item \textbf{Variability and Flexibility:} By adjusting the value of $K$, one can explore a range of transformations, each with its own set of advantages and trade-offs. The flexibility offered by $K$ allows for optimization based on specific requirements, be it maximum data compression or minimal computational overhead.
\end{itemize}

\subsection*{Local Testability and Error Detection}

One of the most promising aspects of the Set Shaping Theory is its inherent capability for local testability. Local testability refers to the ability to efficiently test the integrity of a piece of data without examining the entire data set \cite{goldreich2017} \cite{cheraghchi2005}. In traditional error correction techniques, redundancy is incorporated within the codewords. Set Shaping Theory offers an alternate method, where we introduce redundancy at the message level itself.
\newline\newline
The essence of this method revolves around the bijection functions $f(X^N) = Y^{N+K}$ that transform one set of strings into another of equal size but with strings of greater length. The transformation process conditions the emission probability of the dependent variable $y$ based on previously emitted variables. As a consequence, if a decoder encounters a symbol with a conditional probability of zero, it is an immediate indication of an error in the message, facilitating local testability \cite{testability_kozlov}.
\newline\newline
Several advantages arise from this approach:
\begin{itemize}
    \item \textbf{Efficient Local Testing}: Due to the reshaped strings in $Y^{N+K}$ being chosen based on minimized average information content, local testing becomes inherently more efficient.
    \item \textbf{Reduced Entropy and Error Detection}: While the strings in $Y^{N+K}$ are elongated, their selection is based on reduced entropy. This ensures not only better compression but also enhances the efficiency of local testability.
    \item \textbf{Intrinsic Error Detection}: The conditioning of emission probabilities allows for swift error detection. Encountering a symbol with a conditional probability of zero directly signals a fault in the message, making the method intrinsically robust for error detection.
\end{itemize}

\subsubsection*{Impact of K on testability}
The parameter $K$ in the Set Shaping Theory is pivotal not just for the transformation dynamics but also for enhancing the local testability of data. By studying the variations in average information content with respect to $K$, we can glean insights into the effectiveness of local testability and error detection.
\newline

In an effort to comprehend the practical application of the $f$ function, especially in the design of locally testable codes, a comprehensive study was conducted to investigate how the average information content fluctuates concerning the parameters $|A|$ and $K$. A source $X = (x; A; P)$ with $|A| = 3$ was employed to generate strings of length 100 using a uniform probability distribution. The transformation function $f$ was then applied to these strings with $K$ varying from 1 to 5 \cite{testability_kozlov}. This experiment was conducted by employing the Monte Carlo method \cite{monte_carlo}, simulating a vast set of up to $10^7$ strings to ensure accurate results.
\newline\newline
From the results, a few key observations can be made:

\begin{itemize}
    \item \textbf{Enhanced Error Sensitivity}: As the value of $K$ increased in the Set Shaping Theory's transformation process, the transformed data exhibited heightened sensitivity to errors or anomalies. This means that the reshaped data became more adept at revealing inconsistencies, making the detection of errors or deviations more straightforward and efficient. This heightened sensitivity offers promising implications for systems where swift and precise error detection is paramount \cite{testability_kozlov}.

    \item \textbf{Trade-offs with $K$}: While the benefits of increasing $K$ are evident in terms of enhanced error sensitivity, there's an inherent balancing act. Beyond a certain point, the advantages may start to plateau or even reverse. This suggests that while $K$ can be used to fine-tune the transformation process, it's essential to determine its optimal value to harness the full potential of the Set Shaping Theory without incurring undue computational costs \cite{testability_kozlov}.
\end{itemize}

\subsection*{Overcoming Inefficiencies in Entropy Coding}
One of the foundational premises of the Set Shaping Theory (SST) is its approach to addressing inefficiencies inherent in the coding of individual random sequences when the source isn't known. It's crucial to understand that SST doesn't promise to compress a random sequence. Instead, it refines the coding scheme and the coded sequence, thus optimizing the process and mitigating the inefficiencies present.
\newline\newline
In traditional coding schemes, when encoding an individual random sequence without prior knowledge of the source, inherent inefficiencies creep in. These inefficiencies arise from the lack of context or understanding of the sequence's source, leading to sub-optimal encoding \cite{entropy_inefficiency}.
\newline
SST offers a paradigm shift by reshaping the sequence, guided by the parameter $K$, to address and minimize these inefficiencies. This reshaping doesn't compress the sequence but rather restructures it, making it more amenable to efficient encoding.
\newline
Through the principles of SST, both the coding scheme and the coded sequence undergo optimization. By focusing on sequences with lesser information content (and elongating them through the shaping order $K$), SST ensures that the subsequent encoding process is more efficient and less prone to the traditional pitfalls of sequence encoding without source knowledge.

\section*{Future Directions of Research}

The Set Shaping Theory, although promising, is still in its early stages, and much remains to be explored to realize its full potential. A pivotal aspect that warrants further research is the parameter $K$ itself.

\subsection*{Exploration of Different Values of $K$}
While our current understanding of $K$ has provided valuable insights into the transformation dynamics, a systematic exploration of various values of $K$ is imperative. A step-by-step increment of $K$ can help in understanding the incremental benefits or downsides. This can assist in identifying an optimal range for $K$ for specific applications.

\subsection*{Algorithmic Optimizations}
As the Set Shaping Theory gains traction, there will be a need for more efficient algorithms to implement the transformations. Developing algorithms that can quickly compute the transformations for larger datasets and higher values of $K$ will be crucial.

\subsection*{Real-world Applications and Benchmarks}
Theoretical explorations need to be complemented by real-world applications, testing the Set Shaping Theory in practical data compression scenarios will validate its efficacy. Establishing benchmarks for different $K$ values across varied datasets can provide a guideline for practitioners.

\section*{Conclusion}
The Set Shaping Theory, with the $K$ parameter at its core, presents a novel approach to data transformation and compression. While at first glance adding length to data strings might seem paradoxical, the underlying principle of focusing on strings with reduced information content makes it a potential game-changer in data compression scenarios. Understanding the intricacies of the $K$ parameter is essential to harness the full potential of this theory. With ongoing research and practical applications, the Set Shaping Theory, guided by $K$, promises to revolutionize data compression methodologies, offering optimized, efficient, and adaptable solutions.

\end{document}